\documentclass[reprint,superscriptaddress,showpacs,amsmath,amssymb,aps,pre,floatfix,shortbibliography]{revtex4-1}
\usepackage{amsthm}
\usepackage{amsmath}
\usepackage[utf8]{inputenc}	
\usepackage[english]{babel}
\usepackage{amsmath}
\usepackage{graphicx}		
\usepackage{natbib}
\usepackage{textcomp}
\usepackage{gensymb}
\usepackage[usenames,dvipsnames,svgnames,table]{xcolor}
\usepackage[hidelinks,colorlinks=false,urlcolor=Cerulean,citecolor=black]{hyperref}
\usepackage{siunitx}
\usepackage{mathrsfs}
\usepackage{multirow}
\usepackage{xfrac}
\usepackage{comment}
\usepackage[normalem]{ulem}

\renewcommand{\i}{\mathrm{i}}

\newcommand{\re}{{\text{re}}}
\newcommand{\im}{{\text{im}}}

\begin{document}

\title{Geometry unites synchrony, chimeras, and waves in nonlinear oscillator networks}

\author{Roberto C. Budzinski}
\thanks{These authors contributed equally}
\author{Tung T. Nguyen}
\thanks{These authors contributed equally}
\affiliation{Department of Mathematics, Western University, London, ON N6A 3K7, Canada}
\affiliation{Brain and Mind Institute, Western University, London, ON N6A 3K7, Canada}
\author{Jacqueline Đo\`{a}n}
\affiliation{Department of Mathematics, Western University, London, ON N6A 3K7, Canada}
\affiliation{Brain and Mind Institute, Western University, London, ON N6A 3K7, Canada}
\author{J\'an Min\'a{\v c}}
\affiliation{Department of Mathematics, Western University, London, ON N6A 3K7, Canada}
\author{Terrence J. Sejnowski}\email{terry@salk.edu}
\affiliation{The Salk Institute for Biological Studies, La Jolla, CA 92037, USA}
\affiliation{Division of Biological Sciences, University of California at San Diego, La Jolla, CA 92093, USA}
\author{Lyle E. Muller}\email{lmuller2@uwo.ca}
\affiliation{Department of Mathematics, Western University, London, ON N6A 3K7, Canada}
\affiliation{Brain and Mind Institute, Western University, London, ON N6A 3K7, Canada}

\begin{abstract}
One of the simplest mathematical models in the study of nonlinear systems is the Kuramoto model, which describes synchronization in systems from swarms of insects to superconductors. We have recently found a connection between the original, real-valued nonlinear Kuramoto model and a corresponding complex-valued system that permits describing the system in terms of a linear operator and iterative update rule. We now use this description to investigate three major synchronization phenomena in Kuramoto networks (phase synchronization, chimera states, and traveling waves), not only in terms of steady state solutions but also in terms of transient dynamics and individual simulations. These results provide new mathematical insight into how sophisticated behaviors arise from connection patterns in nonlinear networked systems.
\end{abstract}

\maketitle

\begin{quotation}
\textbf{The collective behavior of nonlinear oscillator networks has been used to study systems ranging from biology to physics. In this context, the Kuramoto Model (KM) is of great importance. However, it remains difficult to directly relate the structure of a specific network adjacency matrix to the resulting dynamics in nonlinear systems. Here, we use a complex-valued matrix formulation for the KM, whose argument has a correspondence with the original model. This allows us to obtain an analytical approach for the transient dynamics in individual simulations of Kuramoto oscillator networks. We then apply this to study phase synchronization, chimera states, and waves. Our approach gives us a new, geometric perspective of synchronization phenomena in terms of complex eigenmodes, which in turn offers a unified geometry for synchrony, chimera states, and waves in nonlinear oscillator networks.}
\end{quotation}

Networks of nonlinear oscillators have attracted interest across physics, neuroscience, biology, and applied mathematics as models where order can arise without a central coordinator \cite{strogatz2001exploring, strogatz2018nonlinear}. Oscillator networks have been applied to study the behavior of insects \cite{buck1988synchronous,ermentrout1991adaptive}, patterns of social behavior \cite{pluchino2006compromise,pluchino2006opinion}, neural systems \cite{breakspear2010generative,cabral2011role}, and physical systems \cite{wiesenfeld1998frequency,valagiannopoulos2020injection}. In this domain, the Kuramoto oscillator has emerged as the central mathematical model for studying these synchronization phenomena \cite{acebron2005kuramoto, kuramoto1975self, crawford1994amplitude, crawford1995scaling, rodrigues2016kuramoto}. While it is quite simple, the Kuramoto model has led to the discovery of new dynamical phenomena well beyond the initial study of nodes evolving to the same phase. Three major dynamical phenomena studied in this model are complete phase synchronization \cite{kuramoto1975self}, partially synchronized ``chimera'' states \cite{kuramoto2002coexistence,abrams2004chimera,abrams2006chimera,hagerstrom2012experimental}, and traveling waves \cite{paullet1994stable}. While much study in mathematics has focused on the existence and stability of these states, the link between network topology (i.e.~the structure of connections) and transient network dynamics remains unclear. Here we report a new mathematical approach to nonlinear oscillator networks that can provide insight into how the structure of an individual network drives synchronization phenomena for any single simulation. We find this new approach can explain the major synchronization phenomena that have been found in the Kuramoto model, not only in terms of stable, steady state solutions but also in terms of transient behaviors in finite time.

The original Kuramoto model (KM) with phase-lag is given by:
\begin{equation}
    \dot{\theta_{i}} = \omega_{i} + \epsilon \sum\limits_{j=1}^{N} A_{ij} \sin{(\theta_{j} - \theta_{i} - \phi)},
    \label{eq:main_kuramoto}
\end{equation}
where $\theta_{i}(t) \in \mathbb{R}$ is the phase of the $i^{\mathrm{th}}$ oscillator at time $t$, $\omega_{i}$ is its natural frequency, $N$ is the number of oscillators, $\epsilon$ scales the coupling strength, $\bm{A}$ is the adjacency matrix, and $\phi$ is the phase lag. We focus here on the case where all oscillators have the same natural frequency ($\omega_i = \omega~\text{for all}~i \in [1,N]$). The value of $\phi$ transforms the interaction function from the standard attractive case ($\phi = 0$) to the cosine version with neutral coupling ($\phi = \sfrac{\pi}{2}$) \cite{gong2019repulsively}.

Following recent work \cite{muller2021algebraic}, we can use a complex-valued, algebraic approach to the Kuramoto model. Starting with Eq. (\ref{eq:main_kuramoto}), we can change to a rotating coordinate frame \cite{strogatz1988collective} and set $\omega = 0$ without loss of generality. By subtracting an additional imaginary component in the interaction term,
\begin{equation}
\dot{\psi}_i=\epsilon \sum_{j=1}^N A_{ij} \big[\sin(\psi_j-\psi_i-\phi)-\i \cos(\psi_j-\psi_i-\phi) \big],
\label{eq:analytical_KM}
\end{equation} 
we arrive at a complex-valued system in a new variable $\psi_i \in \mathbb{C}$ that provides an analytical, algebraic approach to the dynamics of the original nonlinear KM. Specifically, while the addition of the imaginary component results in a modified system, using the analytical solution of Eq.~(\ref{eq:analytical_KM}) (detailed below) to propagate across multiple time windows produces trajectories that match those in the original KM for long timescales (Supplementary Figs. S1, S2, and S3). It is important to note that this specific form of the additional imaginary component is necessary (cf.~Supplemental Materials, Sec.~I.B). 

We now show that Eq.~(\ref{eq:analytical_KM}) has a closed-form, analytical solution. To do this, we can now multiply by $\i$ and use Euler's formula to obtain (cf.~Supplement Sec. I B):
\begin{equation}
\i \dot{\psi}_i= \epsilon e^{-\i \psi_i} \sum_{j=1}^N A_{ij} e^{-\i \phi} e^{\i \psi_j}.
\end{equation}
Letting $x_i = e^{\i \psi_i}$, we can write:
\begin{equation} 
\dot{x}_i = \epsilon \sum_{j=1}^N A_{ij} e^{-\i \phi} x_{j},
\label{eq:analytical_x}
\end{equation}
which has the general solution:
\begin{equation}
\vec{x}(t)= e^{t \bm{K}} \vec{x}(0),
\label{eq:solution_analytical_model}
\end{equation}
where $\vec{x}=(x_1, \cdots, x_N)$ and
\begin{equation}
    K_{ij}=\epsilon e^{-\i \phi} A_{ij}\,,
    \label{eq:matrix_k}
\end{equation}
where $\bm{K}$ now collates the network topology $\bm{A}$, the coupling strength $\epsilon$, and the phase lag $\phi$. 

We now have two distinct dynamical systems, the original KM and the complex-valued system with the explicit solution in Eq.~(5). In the complex-valued system, $\vec{x} \in \mathbb{C}^N$ has complex-valued elements $x_i(t) \in \mathbb{C}$ whose argument we compare with the numerical solution of the original KM $\theta_i(t) \in \mathbb{R}$ (Supplementary Video 1). In more detail, let $\vec{\psi}=\vec{\psi}_{\re}+\i \vec{\psi}_{\im}$ be the decomposition of $\vec{\psi}$ into the real and imaginary parts. Then we have
\begin{equation}
\vec{x}=e^{\i \vec{\psi}_{\re}- \vec{\psi}_{\im}}=e^{-\vec{\psi}_{\im}} e^{\i \vec{\psi}_{\re}}.
\end{equation}
From this, we can see that the argument of $\vec{x}(t)$ is the real part of $\vec{\psi}$. With this, we can show that the real and imaginary parts of the complex-valued model follow the dynamics (see Supplement Sec.~I.B):
\begin{equation} \label{eq:analytical_KM_real_part}
     \frac{d(\psi_i)_{\mathrm{re}}}{dt}= \epsilon \sum_{j=1}^N A_{ij} \frac{|x_j|}{|x_i|} \sin((\psi_j)_{\mathrm{re}}-(\psi_i)_{\mathrm{re}}-\phi),
\end{equation}
and
\begin{equation} \label{eq:analytical_KM_imag_part}
    -\frac{d(\psi_i)_{\mathrm{im}}}{dt}= \epsilon \sum_{j=1}^N A_{ij} \frac{|x_j|}{|x_i|} \cos((\psi_j)_{\mathrm{re}}-(\psi_i)_{\mathrm{re}}-\phi). 
\end{equation}
Note that Eqs.~(\ref{eq:analytical_KM_real_part}) and (\ref{eq:analytical_KM_imag_part}) show that the complex-valued system is not identical to the original KM since moduli $|x_i(t)|$ in the complex-valued system can exhibit amplitude dynamics that increase in the complex plane. Surprisingly, however, when initialized with unit-modulus initial conditions $|x_i(0)| = 1~\text{for all}~i$, with complex arguments $\mathrm{Arg}[x_i(0)]$ that match the initial phases $\theta_i(0)$ in the original KM, the trajectories in the original and complex-valued KM correspond for a non-trivial window of time (Fig. S1). In this work, we introduce an approach to evaluate Eq.~(\ref{eq:solution_analytical_model}) in short time windows, with the goal of utilizing the analytical form to generate insight into the underlying mechanism of synchronization, chimeras, and traveling waves in nonlinear oscillator networks. The approach involves: (i) taking a unit-modulus state vector $\vec{x}(t_0)$ at the beginning of each window, (ii) propagating the analytical expression by evaluating the matrix exponential in Eq.~(\ref{eq:solution_analytical_model}) for a short time window (down to 1 millisecond in the most challenging case of the chimera state), and then (iii) taking the argument of the result ($\mathrm{Arg}[(\vec{x})_i(t)]~\text{for all}~i \in [1,N]$), as done when we compare the solution of the complex-valued system with the numerical solution of the original KM. With this approach, the trajectories in the complex-valued model match those in the original Kuramoto model in a variety of conditions and across realizations over random initial conditions (Figs.~S6, S7, S8, and S9). A detailed comparison of dynamics in these two models is provided in the Supplementary Materials (Sec.~I.B).

We thus compare the argument of our analytical expression (Eq.~(\ref{eq:solution_analytical_model}), hereafter denoted ``analytical'') to computer simulations of the original nonlinear Kuramoto model (Eq.~(\ref{eq:main_kuramoto}), ``original KM'') throughout the rest of this work. The numerical simulations were performed using Euler's method, and no significant differences were observed using a different integration technique. Importantly, we note that initial conditions (at $t =0$) in the complex-valued Kuramoto model always have unit modulus $|x_i(0)| = 1~\text{for all}~i$, with arguments $\mathrm{Arg}[x_i(0)]$ equal to the initial angles $\theta_i(0)$ in the original KM.

In general (see Supplementary Material - Sec. I.B - for specific conditions), we can write the temporal dynamics of $\vec{x}$ in terms of the eigenvectors and eigenvalues of $\bm{K}$:
\begin{equation}
    \vec{x}(t) = \underbrace{c_{1}e^{\lambda_{1}t}}_{\mu_{1}}\vec{v}_{1} + \underbrace{c_{2}e^{\lambda_{2}t}}_{\mu_{2}}\vec{v}_{2} + \cdots + \underbrace{c_{N}e^{\lambda_{N}t}}_{\mu_{N}}\vec{v}_{N}, 
    \label{eq:solution_analytical}
\end{equation}
where $\lambda_{k}$ is the eigenvalue associated with the eigenvector $\vec{v}_{k}$. The projection on to the $k^{\mathrm{th}}$ eigenmode is given by the inner product $\mu_{k}(t) = \langle \vec{x}(t), \vec{v}_{k} \rangle$. 

This expression now allows us to understand the three major synchronization phenomena previously discovered in the Kuramoto model from a comprehensive mathematical perspective, in terms of the geometry of the $\mu_k$ representing the contribution of the $k^{\mathrm{th}}$ eigenmode. Furthermore, when $\bm{A}$ follows the same connectivity rule for all nodes (i.e.~where connections can be fully specified by a single vector cyclically permuted across all nodes, resulting in a circulant matrix), the circulant diagonalization theorem (CDT) allows obtaining the eigenspectrum analytically, with eigenvectors following an ordering of Fourier modes from low to high spatial frequencies (Supplementary Materials, Sec.~I.B, Figs.~S10 and S13) \cite{davis1979}. Using Eq.~(\ref{eq:solution_analytical}), we can then analyze the system in terms of the eigenmode contributions $\mu_{k}$.

\begin{figure}[t]
    \centering
    \includegraphics[width=0.98\columnwidth]{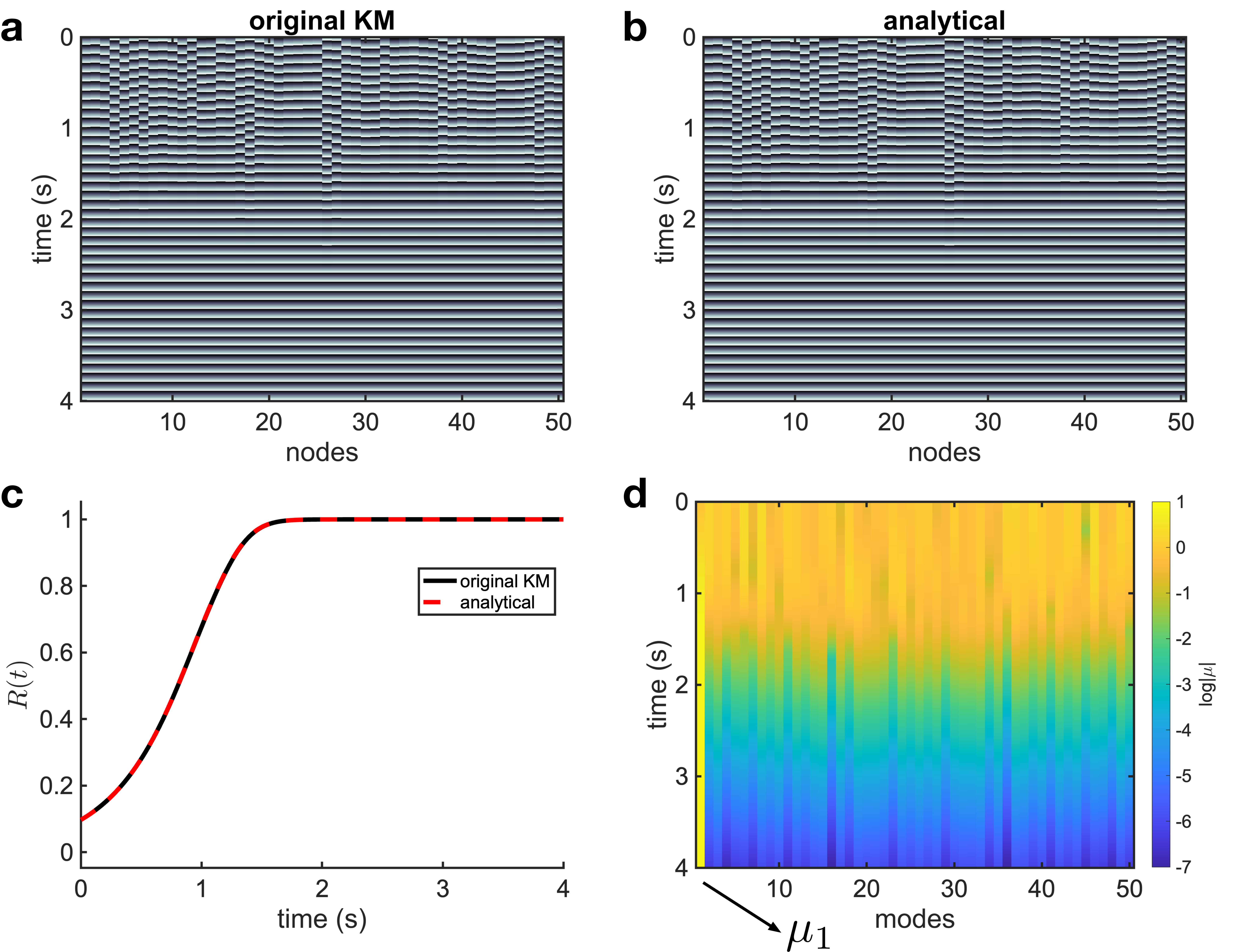}
    \caption{Analytical description of the transition to synchrony. Solutions of the original Kuramoto model (numerical simulation) (a) and the analytical evaluation (b) are plotted in color-code (dark tones indicate values close to $-\pi$ and light tones to $\pi$). As time evolves, the spatiotemporal dynamics become coherent and phase synchronization appears. A quantitative comparison between the two models is provided in Fig.~S1. The level of synchronization, measured by the Kuramoto order parameter (c), starts at a low value with the random initial conditions and quickly approaches unity, which indicates the phase synchronized state of the system. The eigenmode contribution ($\log{|\mu|}$) is plotted in color (d). The contribution of the first eigenmode, which represents zero phase difference across oscillators, dominates when the phase synchronized state appears.}
    \label{fig:synchronization}
\end{figure}

We first studied the emergence of phase synchronized states in Kuramoto oscillators with attractive coupling ($\phi = 0$), where the adjacency matrix $\bm{A}$ is given by a complete graph on $N = 50$ nodes. Figure \ref{fig:synchronization} shows the precise correspondence between (a) the original Kuramoto and (b) the analytical evaluation during the transition from random initial conditions ($\theta_{i}(0) = \mathrm{Arg}[x_i(0)] \sim \mathcal{U}[-\pi,\pi]$) to synchrony. This transition is measured by the Kuramoto order parameter $R(t) = N^{-1} |\sum_{j=1}^N \exp[\i\theta_j(t)]|$ (where the phases are obtained from the original KM or complex-valued model), which also exhibits a precise correspondence between the original KM and analytical evaluation (Fig.~\ref{fig:synchronization}c). Using this analytical approach, we can then analyze the system in terms of the contribution of individual eigenmodes, represented in color in Fig.~\ref{fig:synchronization}d ($\log{|\mu|}$). Prior to the transition to synchrony, the eigenmode contributions are almost uniform. After the transition, however, the contributions shift away from uniformity, and the first eigenmode becomes dominant (yellow line at $\mu_{1}$).

The analytical approach we have introduced provides an opportunity to understand this transition geometrically. Specifically, we can study how the first eigenmode, which represents the synchronized state, behaves in relation to the other modes, which represent more complicated configurations of phase. As the system transitions to synchronization, the contribution of the first eigenmode increases in magnitude (modulus) in the complex plane, while the other eigenmode contributions collapse to the origin (Supplementary Video 2). Because it is the argument of Eq.~(\ref{eq:solution_analytical_model}) that determines the phase dynamics in the analytical approach, the contribution of the first eigenmode, which is associated with the eigenvector $\vec{v}_{1} = (1, 1, \cdots, 1)^{T}$ (and $\mathrm{Arg}[(\vec{v}_{1})_{i}] = 0~\text{for all}~i \in [1, N]$), explicitly brings the network towards the synchronized state. This analysis provides a novel geometric insight into the transition to synchrony and how the pattern of connections in a network determines the transition from an incoherent to a coherent state during an individual simulation. Further, when $\phi = \sfrac{\pi}{2}$, interactions in the network are no longer attractive, and in this case, Kuramoto networks do not synchronize \cite{kuramoto1984cooperative}. Consistent with this, our analysis shows the eigenmode contributions remain uniform across time (Fig.~S5 and Supplementary Video 3). These results demonstrate the geometric view provided by the eigenmodes can provide insight into the transient dynamics of synchronized and incoherent states. 

\begin{figure}[t]
    \centering
    \includegraphics[width=1.0\columnwidth]{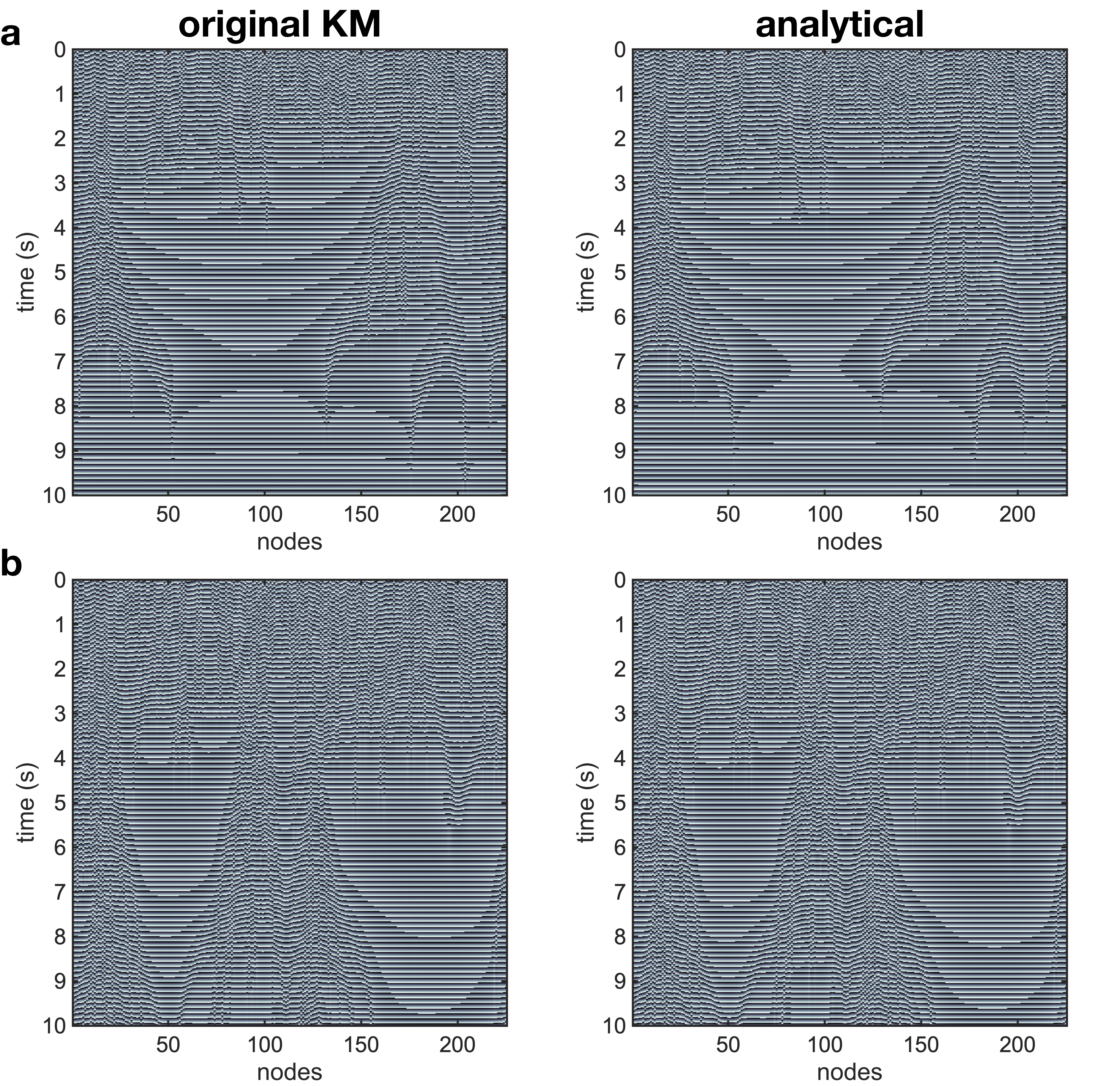}
    \caption{Analytical description of the chimera state. Spatiotemporal dynamics for original KM (simulation) and analytical evaluation are plotted in color-code for the distance-dependent networks with $\phi = 1.15$ (a), and $\phi = 1.30$ (b). In these cases, chimera states are observed, where part of the network is synchronized and coexists with an incoherent region. Here, it is important to emphasize that the original KM and analytical expression are evaluated separately, using only the same initial condition. This demonstrates that the analytical approach is able to capture important details of the Kuramoto dynamics.}
    \label{fig:chimeras}
\end{figure}

We now use this approach to understand two more sophisticated dynamical phenomena that have been discovered in the Kuramoto model:~partially synchronized chimera states and traveling waves. Chimera states, which are a sophisticated mixture of synchronized and non-synchronized clusters in oscillator networks, are known to arise in models with distance-dependent (non-local) connections and a constant phase lag ($\phi$) \cite{parastesh2020chimeras, abrams2004chimera, shanahan2010metastable, wolfrum2011chimera}. Here we consider the case where connections follow a deterministic, distance-dependent power rule that specifies a real-valued connection weight (Supplementary Material, Sec.~I.E.2). Figure \ref{fig:chimeras} depicts the spatiotemporal dynamics for the original KM and our analytical solution for the cases $\phi = 1.15$ (a) and $\phi = 1.30$ (b). In  Fig.~\ref{fig:chimeras}a, one can see a transient chimera, where the system transitions from incoherence, to a chimera state, and finally to phase synchronization. In Fig.~\ref{fig:chimeras}b, the system transitions from incoherence to a chimera state that continues for all times we simulated. Importantly, in both cases the the spatiotemporal dynamics produced by our analytical evaluation depicts a good correspondence with the chimera states observed in the numerical simulation (Supplementary Videos 4 and 5). A quantitative comparison between the analytical and original KM results in this case is provided in Fig.~S3.

\begin{figure*}[htb]
     \centering
     \includegraphics[width=0.85\textwidth]{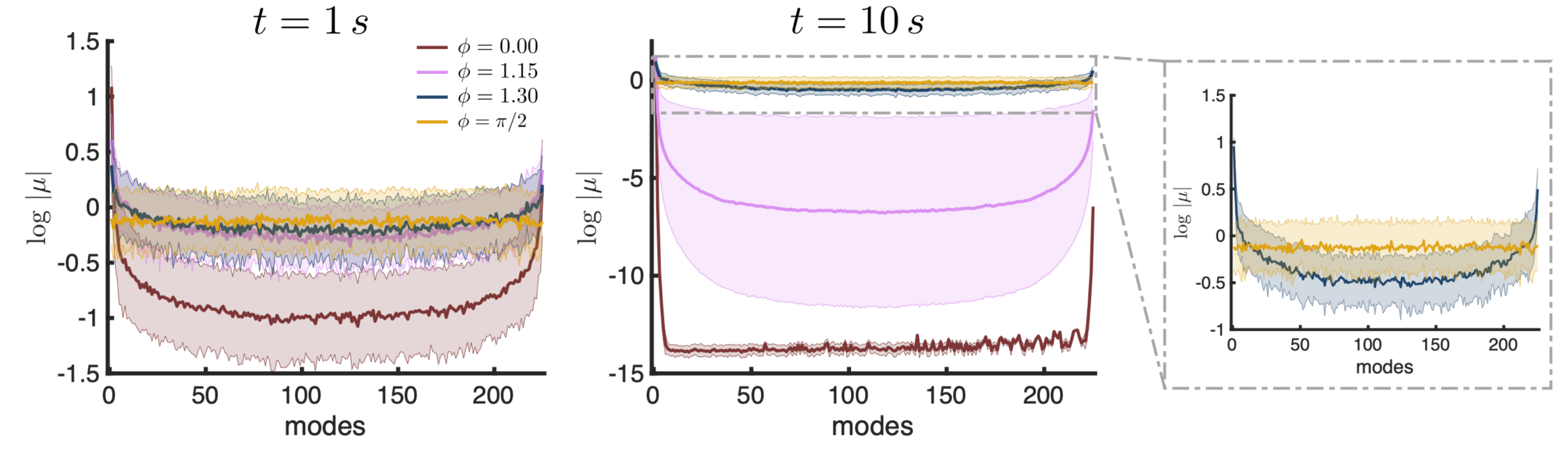}
     \caption{Geometry of eigenmodes explains the mechanism for the chimera state. The eigenmode contribution ($\log{|\mu|}$) for increasing values of $\phi$ at $t = 1$ s (left) and $t = 10$ s (middle) demonstrates the evolution of the first modes in the case of $\phi = 0$, while for $\phi = 1.30$, the first and last eigenmodes contribute in balance (right). Solid lines represent average eigenmode contribution over 100 simulations, and shaded regions represent standard deviation. }
     \label{fig:eigenmodes_chimera}
\end{figure*}

We now use the geometric approach to the Kuramoto dynamics to derive insight into the mechanism underlying chimera states in oscillator networks. Specifically, in networks where the same connectivity rule is applied to each node, the eigenmodes higher than $\mu_{1}$ take the form of traveling waves from low to high spatial scales in their arguments (Fig.~S13) \cite{davis1979}. Analysis of the eigenmode contributions ($\log{|\mu|}$), plotted as solid lines representing an average over $100$ simulations (with shaded regions representing the standard deviation over simulations), reveals the chimera is created by an interplay of the synchronizing mode $\mu_{1}$ and a set of modes representing waves traveling in opposite directions with progressively higher spatial frequency ($t = 10$ s, blue line, Fig. \ref{fig:eigenmodes_chimera}, right panel). This is in contrast to the completely synchronous state, where the contribution of the first eigenmode and the rest differ by more than $12$ orders of magnitude $10$ seconds into the simulation (burgundy line, Fig. \ref{fig:eigenmodes_chimera}, middle panel), and in contrast to the incoherent case ($\phi = \sfrac{\pi}{2}$), where the contributions across eigenmodes remain uniform (yellow line, Fig. \ref{fig:eigenmodes_chimera}, right panel). These results allow us to understand the emergence of the chimera as an interplay between specific types of modes in the Kuramoto system. Furthermore, analysis of the aggregate connectivity matrix $\bm{K}$ illustrates the underlying mechanism for this interplay of modes is a rotation of the eigenvalues in the complex plane, which reduces the difference between the real part of the eigenvalue associated with the synchronized state and the real part of the rest (Figs.~S10 and S11).

\begin{figure}[htb]
    \centering
    \includegraphics[width=0.95\columnwidth]{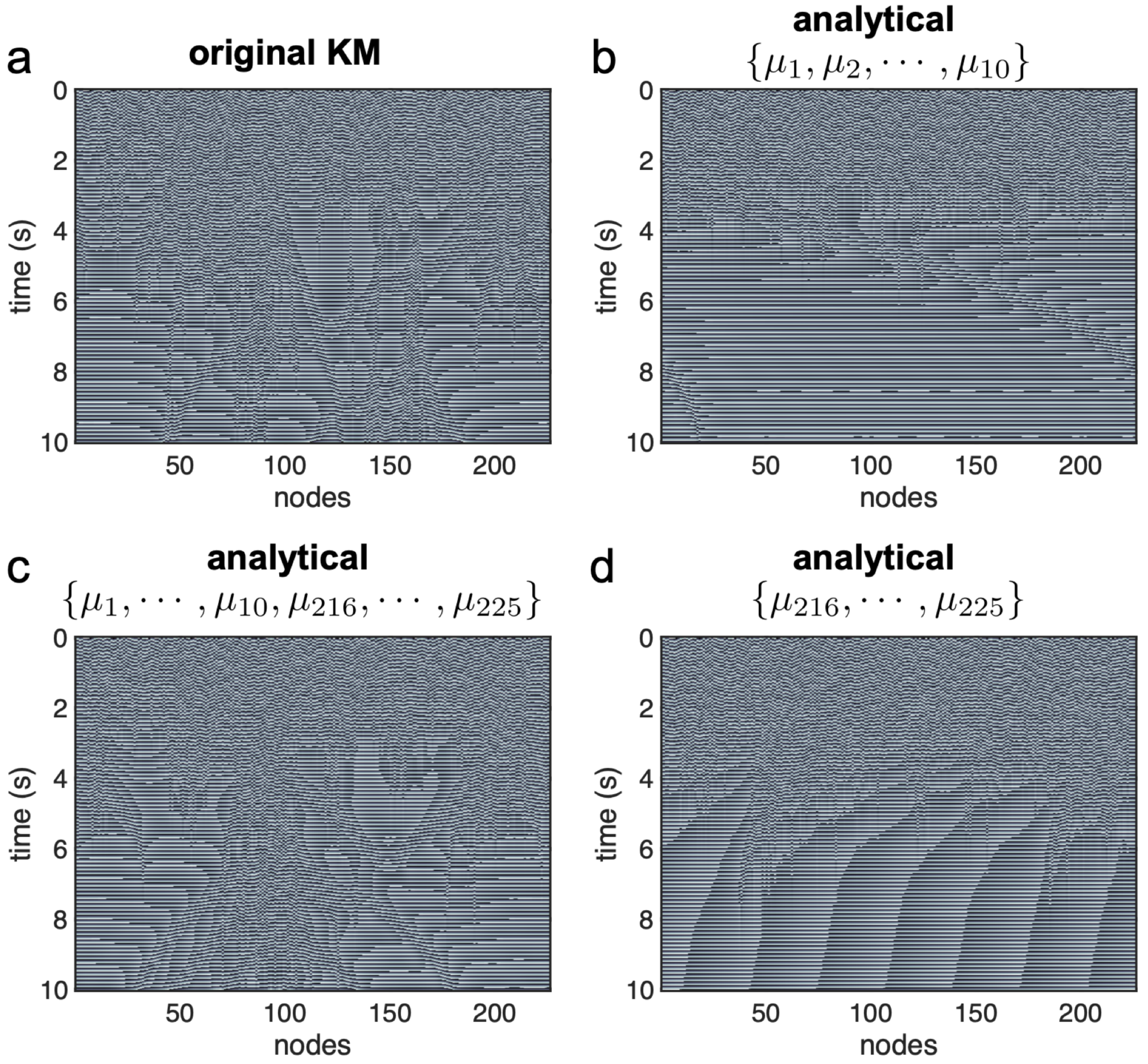}
    \caption{A subset of contributing eigenmodes creates an analytical approximation to the chimera state numerically simulated in the original KM (a). The first 10 eigenmodes $\{ \mu_{1}, \mu_{2}, \cdots, \mu_{10} \}$ create a partially synchronized region that is too broad (b), while the last 10 eigenmodes $\{ \mu_{216}, \mu_{217}, \cdots, \mu_{225} \}$ create no signature of a chimera (d). However, using the first 10 and last 10 eigenmodes $\{\mu_{1}, \mu_{2}, \cdots, \mu_{10}, \mu_{216}, \mu_{217}, \cdots, \mu_{225} \}$ creates an analytical approximation to the chimera state (c).}
    \label{fig:approximation_chimera}
\end{figure}

To illustrate the insight into the chimera state this geometric approach can provide, we can now rewrite the analytical solution using a subset of contributing modes and compare this truncated approximation to the numerical simulation (Fig.~\ref{fig:approximation_chimera}a, see also Fig.~S12). With only the first $10$ modes in the solution ($\{ \mu_{1}, \mu_{2}, \cdots, \mu_{10} \}$), the synchronization in the system is overestimated (Fig. \ref{fig:approximation_chimera}b). Moreover, with only the last $10$ modes in the solution ($\{ \mu_{216}, \mu_{217}, \cdots, \mu_{225} \}$), the synchronization is underestimated (Fig.~\ref{fig:approximation_chimera}d), leaving the system with no signature of the chimera. With the first and last 10 modes ($\{\mu_{1}, \mu_{2}, \cdots, \mu_{10}, \mu_{216}, \mu_{217}, \cdots, \mu_{225} \}$), however, the dynamics in the numerical simulation are recovered, and the main structures defining the chimera state are observed in the analytical approximation (Fig. \ref{fig:approximation_chimera}c). These results demonstrate this analytical approach can capture highly non-trivial dynamical phenomena such as chimera states that emerge from the network structure and nonlinear dynamics of the Kuramoto model.

\begin{figure}[bth]
    \centering
    \includegraphics[width=0.95\columnwidth]{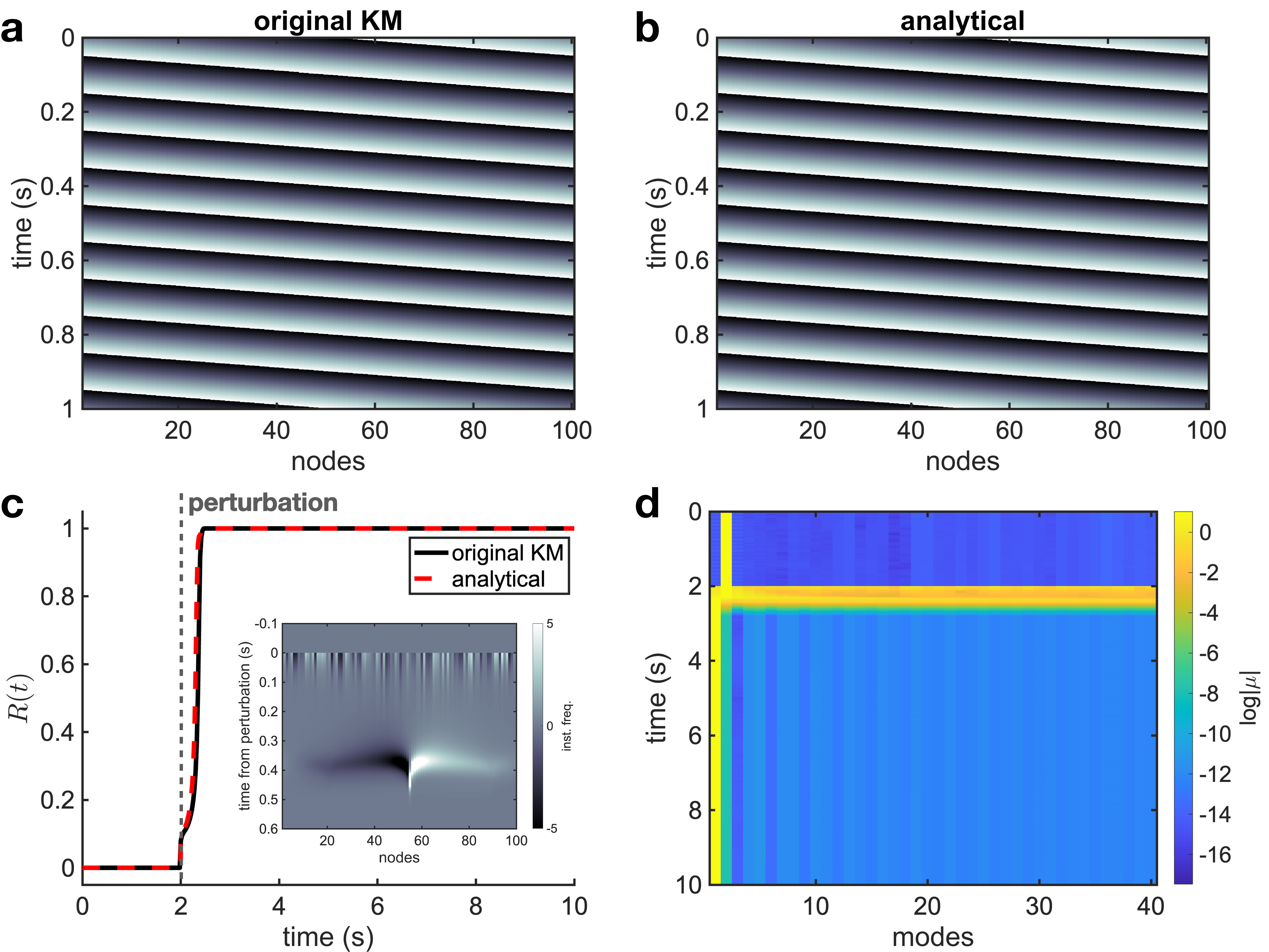}
    \caption{Original KM (a) and analytical evaluation (b) for an example of a wave state on a ring graph ($N = 100, k = 1$, see Supplement Sec. I E 2), resulting from a special set of initial conditions, a “twisted” state (see Eq. (27) in the Supplement). In this context, traveling wave states constitute the contribution of a single eigenmode to the dynamics ($2^{\mathrm{nd}}$ eigenmode, in this case) (cf.~Fig. S13). The Kuramoto order parameter (c) in the traveling wave state remains at $R = 0$. However, when a finite perturbation is applied to the system at $t = 2$ s, the network transitions to a phase synchronized state ($R = 1$). The relative instantaneous frequencies following the perturbation exhibit a non-trivial self-organized state during the transition to synchrony (c, inset). This transition is captured by the eigenmode contributions (d), where the $2^{\mathrm{nd}}$ eigenmode gives way to the $1^{\mathrm{st}}$ eigenmode when phase synchronization is reached.}
    \label{fig:traveling_waves}
\end{figure}

Finally, we considered traveling waves in the Kuramoto model. While, when $\phi = 0$, the synchronized state occurs starting from a broad range of random initial conditions (as in Fig. \ref{fig:synchronization}), Kuramoto networks can also exhibit traveling wave states that exist for arbitrary lengths of time (original KM and analytical, Fig. \ref{fig:traveling_waves}a and Fig. \ref{fig:traveling_waves}b. A quantitative comparison between these two results is provided in Fig.~S4). While the stability of such traveling wave, or ``twisted'', states has been the subject of much investigation \cite{kopell1986symmetry,taylor2012there,deville2016phase,townsend2020dense}, our geometric approach provides insight into transient wave dynamics, which have recently been observed to play important roles in several fields, from neural systems \cite{ermentrout2001traveling,muller2018cortical} to ecology \cite{morozov2020long} and power grids \cite{halekotte2021transient}. Analysis of eigenmode contributions during a traveling wave on a ring graph reveals these states exist as contributions from a single eigenmode representing a more structured configuration of phase than $\mu_{1}$ (Figs.~S14 and S15).

If the network receives a perturbation with sufficient magnitude, however, the traveling wave state can transition to synchrony (Fig.~\ref{fig:traveling_waves}c). In this case, our analytical approach also captures the timecourse of these transient dynamics well (see Supplementary Video 6 and Fig. S15). In this case, the transition to synchrony has a specific signature: while the eigenmode contribution is localized to $\mu_{2}$ prior to the perturbation, the contributions become spread out during the transition before collapsing to $\mu_{1}$ (Fig.~\ref{fig:traveling_waves}d), since the system transitions from a wave state to phase synchronization (where $R = 1$ in Fig.~\ref{fig:traveling_waves}c). Furthermore, after the perturbation is applied, the instantaneous frequency is no longer the same across oscillators, which results in a non-trivial, self-organized ``falcon” shape in the instantaneous frequencies leading to phase-synchronization after the transition (inset, Fig. \ref{fig:traveling_waves}c). These results demonstrate the geometric view can provide novel insight into transient dynamics in Kuramoto systems following perturbation.

The connection between network structure and resulting nonlinear dynamics is a central question in modern physics. In this work, we have studied an analytical approach to the Kuramoto model, a canonical model for synchronization dynamics in nature. The geometric view enabled by our analytical approach provides insight into three major synchronization phenomena in the Kuramoto model (phase synchronization, partially synchronized chimera states, and traveling waves). The key novel feature of our analytical approach is that it is valid at finite scales and for individual realizations of the model, which can provide more detailed insight into specific, moment-by-moment dynamics in the system than statistical or asymptotic theoretical approaches. Here, the insight provided by this approach allows us to explain the fully synchronized state as a dominant first eigenmode of the complex system (Fig.~\ref{fig:synchronization}), partially synchronized chimera states as an interplay between modes representing a set of waves traveling in opposite directions (Fig.~\ref{fig:chimeras}, Fig.~\ref{fig:eigenmodes_chimera}, and Fig.~\ref{fig:approximation_chimera}), and traveling waves as localizations in higher eigenmodes (Fig.~\ref{fig:traveling_waves}). This unifying insight into three main dynamical phenomena that have been discovered in the Kuramoto model demonstrates the utility of the non-asymptotic, algebraic approach to network structure and nonlinear dynamics in this work. 

Our approach involves a complex-valued system that admits an explicit analytical solution and whose trajectories correspond to the original, nonlinear KM for a non-trivial window of time. Motivated to see how precisely the trajectories between the two models could match, here we developed an approach to evaluate the analytical expression in Eq.~(5) in short time windows. Again, briefly, we start with unit-modulus complex-valued initial conditions, corresponding to the initial phases of the original KM, at the start of each window. We then evaluate the matrix exponential for a short time window (down to 1 millisecond in the most challenging case of the chimera state). Finally, we take the argument of each element in the resulting state vector ($\mathrm{Arg}[(\vec{x})_i(t)]~\text{for all}~i \in [1,N]$), as always done when comparing the solution of the complex-valued system with the numerical solution of the original KM.

This approach to evaluating the complex-valued system represents an initial step towards an operator-based approach to nonlinear dynamics. Specifically, while it is possible to evaluate the complex-valued system directly, by applying the matrix exponential to produce a solution that agrees with the original KM for a short, but non-trivial, length of time, here we find  this windowed approach produces trajectories that precisely match those in the original, nonlinear KM across many different dynamical regimes, including synchronization, chimeras, and traveling waves. The approach taken in this work can be viewed as the combination of two operators:~the matrix exponential and $\mathrm{Arg}$. Importantly, while this approach is more complicated than a standard evaluation of the solution to a linear system by using the matrix exponential, this approach allows us to describe the microscopic evolution of the Kuramoto system in terms of a linear operator acting on an instantaneous state. This, in turn, permits analytical insight into the mechanisms underlying states such as chimeras in terms of the spectrum of the adjacency matrix. This approach thus allows us to connect the nonlinear dynamics generated in an individual simulation of a network of Kuramoto oscillators to the specific adjacency matrix for that system. In addition to investigating the mechanisms of these dynamical phenomena further using spectral graph theory, we aim to understand more fully how the nonlinear dynamics in the original KM can be described to the high precision necessary to match a chimera state, which is known to be a chaotic transient \cite{wolfrum2011chimera}, by the iterative application of two simple operators.

Finally, because the Kuramoto model has been extensively studied both as a model for neural dynamics \cite{mirollo1990synchronization} and as a fundamental model for computation \cite{ermentrout2001traveling,cook1989mean,arenas1994phase}, these results open up several possibilities for studying the connections between network structure, nonlinear dynamics, and computation. The importance of recurrent network structure is becoming increasingly recognized both in biological \cite{kietzmann2019recurrence} and artificial \cite{liang2015recurrent} visual processing, and in recent work, we introduced a theoretical framework for studying how recurrent connections and traveling waves shape computation in the visual system \cite{muller2018cortical}. Understanding how precise network structure can support specific computations in the brain and artificial learning systems through these analytical approaches may lead to a more comprehensive mathematical understanding of neural computation in future work.

\section*{Supplementary Material}

The supplementary material contains additional mathematical details about the methodology, models and networks, complementary figures, description of the parameters used in the analyses, and description of the supplementary movies.

\begin{acknowledgments}
This work was supported by BrainsCAN at Western University through the Canada First Research Excellence Fund (CFREF), the NSF through a NeuroNex award (\#2015276), the Natural Sciences and Engineering Research Council of Canada (NSERC) grant R0370A01, ONR N00014-16-1-2829, NIH EB009282, NIH EB026899, the Swartz Foundation, and SPIRITS 2020 of Kyoto University. J.M.~gratefully acknowledges the Western University Faculty of Science Distinguished Professorship in 2020-2021. The authors would like to thank Dr. Robin Delabays for insightful discussions on this work. In addition, the authors would like to thank both anonymous Referees for their work and comments to improve the manuscript.
\end{acknowledgments}

\section*{Data Availability}
The code and data that support the findings of this study are openly available in \href{http://mullerlab.github.io}{\textcolor{Cerulean}{http://mullerlab.github.io}}.


%

\end{document}